\documentclass[a4paper,10pt]{article}
\usepackage{fullpage}
\usepackage{amsfonts} 
\usepackage{amssymb} 
\usepackage{latexsym} 
\usepackage{epsfig} 
\usepackage{amsmath}
\usepackage{amsthm}
\usepackage[all]{xy}
\usepackage{mathrsfs}
\newcommand{\cc}{\mathbb{C}} 
\newcommand{\oo}{\mathcal{O}} 
\newcommand{\gx}{\mathfrak{X}} 
\newcommand{\gl}{\mathcal{L}} 
\newcommand{\mx}{\mathfrak{m}} 
\newcommand{\lra}{\longrightarrow}%
\newcommand{\ra}{\rightarrow}%
\newcommand{\xra}{\xrightarrow}%

\newcommand{\cvf}[2][]{\frac{\partial #1}{\partial #2}} 
\newtheorem{thm}{Theorem}[section]%
\newtheorem{lem}[thm]{Lemma}%
\newtheorem{cor}[thm]{Corollary}%
\newtheorem{prop}[thm]{Proposition}%
\theoremstyle{definition} 
\newtheorem{defn}[thm]{Definition}%
\newtheorem{exmp}[thm]{Example}%
\newtheorem{rem}[thm]{Remark}%
\DeclareMathOperator{\Lie}{Lie}%
\DeclareMathOperator{\rank}{rk}%
\DeclareMathOperator{\Res}{Res}%

\title{Deformations of functions and F-manifolds}%
\author{Ignacio de Gregorio\footnote{Partially suppported by EPSRC
    grant (00801853), EU mobility project OMATS (HPMT-CT-2000-00104)
    and University of Warwick Research Fellowship Scheme.}}%
\begin{document}
\maketitle%
\abstract{We study deformations of functions on isolated
  singularities. Give a unified proof of the equality of Milnor and
  Tjurina numbers for functions on isolated complete intersections
  singularities and space curves. We define the structure of
  $F$-manifold on the base space of the miniversal deformation in both
  of the above cases. As a corollary, we proved a conjecture of
  V.~Goryunov stating that the critical values of the miniversal
  unfolding of a function on a space curve are generically local
  coordinates on the base space of the deformation.}

\section{Introduction}
\label{sec:introduction}
The theory of Frobenius manifolds plays a central role in Mirror
Symmetry, after the construction by Givental and
Barannikov~({\cite{Bar00}}) of an isomorphism between the small
quantum cohomology of $\mathbb{CP}^n$ and the base space of the
miniversal deformation of the linear function $f=x_1+\dots+x_{n+1}$ on
the divisor $D\colon=\{x_1\cdots{}x_{n+1}=1\}$. There are now a number
of conjectures stating similar isomorphisms between quantum cohomology
rings of algebraic varieties and unfoldings of functions on affine
varieties. In this article we propose a Singularity Theory framework
in which at least one of the ingredients making up the definition of
Frobenius manifolds, namely the multiplication, can be naturally
defined. This structure is known as \emph{$F$-manifold}
(\cite{HerMan99}, \cite{Her02}).

A seemingly inescapable feature of this construction is that the
multiplication is not defined on the whole tangent sheaf of the base
space but only on a certain subsheaf, that of logarithmic vector
fields to the discriminant. Contrary to those Frobenius manifolds
constructed from unfoldings of isolated hypersurface singularities,
our construction does contain some promising candidates to mirrors of
algebraic varieties.

The main result of this article can be stated as follows:
\begin{thm}
  Let $f\colon(X,x)\ra(\cc,0)$ be a function-germ with an isolated
  singularity on a isolated complete intersection or a space curve.
  Then the sheaf $\Theta(-\log\Delta)$ of logarithmic vector fields of
  the discriminant of its miniversal deformation is in a natural
  manner an (logarithmic) $F$-manifold. Moreover, each stratum of the
  logarithmic stratification of the base space inherits this
  structure.
\end{thm}

The content of the article is as follows. First we provide a
construction of the miniversal deformation of a function on a singular
variety. We define a morphism closely related to the Kodaira-Spencer
map that will be used to define the multiplication. Secondly we state
a condition that ensures the equality of the dimension of the
miniversal base space space (Tjurina number) and the number of
critical points of an unfolding of $f$ in the smooth fiber of the
deformation. We then show that the condition holds for functions on
isolated complete intersection singularities and (reduced) space
curves. This provides a unified treatment to the $\mu=\tau$-type
results of V.~Goryunov~(\cite{Gor95}) in the case of functions on
isolated complete intersections singularities and of D.~van~Straten
and the second author in the case of functions on space curves
(\cite{MonStr01}). Our methods are closer to those of \cite{MonStr01}.
To finish, we prove that the multiplication satisfies an integrability
condition, making it into a logarithmic $F$-manifold.

Before going into the technical details we would like to work out a
relatively simple example in which a full Frobenius structure can be
constructed, namely that of function $f=x^p+y^q$ on the ordinary
double point $X\colon{}xy=0\hookrightarrow{}\cc^2$. This case is
closely related to the construction of Frobenius manifolds on Hurwitz
spaces by B.  Dubrovin (\cite{Dub96}, \cite{Ros98}), although as we
are also collapsing the curve a new structure on the discriminant is
made apparent. The aim of this example is firstly to guide the reader
through the rest of the article and secondly to show how indeed our
construction contains some interesting examples in Mirror Symmetry. It
appears to be known among specialists that the resulting Frobenius
manifold is the mirror of the orbifold $\mathbb{CP}(p,q)$.

\subsubsection*{Functions on the double point}

Let us consider a function germ $f=x^p+y^q$ on the $A_1$-singularity
$X\colon{}xy=0$. The result of the calculation that we are going to
carry out can be resumed in the following theorem. We remark that
certain aspects of the proof, particularly the multiplication, will
only be evident after applying the results given in the main body of
this article. The divisor $\Delta$ denotes the discriminant of the
miniversal deformation of $f$.
\begin{thm}
  The sheaf $\Theta(-\log\Delta)$ is naturally endowed with a
  multiplication $\star$ and a flat bilinear pairing $\left<,\right>$
  satisfying
  \begin{displaymath}
    \left<u\star{}v,w\right>=\left<u,v\star{}w\right>~\text{for any $u,v,w\in{}\Theta(-\log\Delta)$}
  \end{displaymath}
  There exists a conformal Euler vector field $E$, that is:
  \begin{displaymath}
      \Lie_E(\star)=\star\text{ and }\Lie_E(\left<,\right>)=\left<,\right>
  \end{displaymath}
\end{thm}
\begin{proof}
  Let us begin by constructing the multiplication. The miniversal
  deformation of $f$ is given by the function
  $F=c+\sum_{i=1}^{p-1}a_{i}x^i+x^p+\sum_{i=1}^{q-1}b_iy^i+y^q$ on the
  fibration $\pi(x,y,a,b,c)=(xy,a,b,c)$, where
  $a=(a_{p-1},\dots,a_{1})$ and $b=(b_{q-1},\dots,b_1)$ (see
  Cor.~\ref{versality-criterion} and the paragraph below). Let
  $\Delta\colon\epsilon=0$ be the (smooth) discriminant of $\pi$. The
  multiplication is defined by the following lifting process: lift a
  vector field $u\in\Theta(-\log\Delta)$ to $\tilde{u}$ such that
  $t\pi(\tilde{u})=u\circ{}\pi$. Differentiating $F$ with respect to
  $\tilde{u}$ we obtain an element in the ring of germs
  $\oo_{\gx,0}=\cc\{x,y,a,b,c\}$. We denote by $t'F(u)$ its class in
  the quotient $\oo_{\gx,0}/(H)$, where $H$ is the Jacobian
  determinant
  \begin{displaymath}
    \cvf[(F,\pi_1=xy)]{(x,y)}=\sum_{i=1}^{p-1}ia_{i}x^{i}+px^p-\sum_{i=1}^{q-1}ib_{i}y^{i}-qy^q
  \end{displaymath}
  It will clear from later constructions (although it can be checked
  directly) that the map $t'F$ so constructed is an isomorphism of
  free modules over the base of rank $p+q$. We use it to pull back the
  algebra structure on $\oo_{\gx,0}/(H)$ so defining a multiplication
  $\star$ in $\Theta(-\log\Delta)$.
  
  To define the multiplication, we consider the relative dualising
  form $\alpha=dx\wedge{}dy/d\pi_1$ and use it to identify
  $\oo_{\gx,0}$ with $\omega_{\gx/B,0}$. Hence we have $dF=H\alpha$
  and consider the Grothendiek residue pairing on
  $\omega_{\gx/B,0}/\oo_{\gx,0}\left<dF\right>$. We use $t'F$ to
  define a multiplicatively invariant non-degenerate bilinear pairing
  on $\Theta(-\log\Delta)$. For $u,v\in\Theta(-\log\Delta)_b$, it is
  explicitly given by
  \begin{displaymath}
    \left<u,v\right>=\int_{\partial{}X_b}\frac{t'F(u)t'F(v)}{H}\alpha
  \end{displaymath}
  \noindent{}being $\partial{}X_b$ the boundary of an appropriate
  representative of the fiber $\pi^{-1}(b)$. For
  $b\in{}B\setminus\Delta$, the fiber $X_b$ is a smooth rational curve
  with two points deleted, say $\infty_1$ and $\infty_2$,
  corresponding to $x=0$ and $y=0$. Hence the pairing can be expressed
  as
  \begin{equation}
    \label{eq:1}
    \left<u,v\right>=-\Res_{\infty_1}\frac{t'F(u)t'F(v)}{H}\alpha+\Res_{\infty_2}\frac{t'F(u)t'F(v)}{H}\alpha
  \end{equation}
  If we take the free basis of $\Theta(-\log\Delta)$ given by
  $\epsilon\cvf{\epsilon}$ and the rest of coordinate vector fields,
  the decomposition (\ref{eq:1}) allows us to express the matrix of
  $\left<,\right>$ as a sum, each summand corresponding to the
  residues at each point. A direct calculation, necessary for what
  follows, shows that the matrix is given by
  \begin{equation}
    \label{eq:2}
    \begin{pmatrix}
      0&0&0&4^{-1}\\
      0&0&0&0\\
      0&0&M_{\infty_1}&0\\
      4^{-1}&0&0&0
    \end{pmatrix}
    +
    \begin{pmatrix}
      0&0&0&4^{-1}\\
      0&M_{\infty_2}&0&0\\
      0&0&0&0\\
      4^{-1}&0&0&0
    \end{pmatrix}
  \end{equation}
  \noindent{}where $M_{\infty_1}=\left(
    \begin{smallmatrix}
      2b_2&3b_3&4b_4&\dots&(q-1)b_{q-1}&q\\
      3b_3&4b_4&\dots&\dots&q&0\\
      4b_4&\dots&\dots&\dots&\dots&0\\
      \dots&\dots&\dots&\dots&\dots&\dots\\
      q&\dots\dots&\dots&\dots&\dots&0\\
    \end{smallmatrix}\right)^{-1}$
  and analogously for $M_{\infty_2}$.
  
  To show that the pairing is indeed flat, we compute flat
  coordinates. Let $b=(\epsilon_0,a_0,b_0)\in{}B\setminus{}\Delta$.
  As $F$ has a pole of order $p$ at $\infty_2$, we can find a local
  coordinate $u$ at $\infty_2$ such that $F=u^{-p}$. On the other
  hand, the function $xu$ is holomorphic and not vanishing at
  $\infty_2$. Fixing a branch of $\log$ we can expand it as a power
  series:
  \begin{displaymath}
    \log{}xu=t_0+t_1u+\dots+t_{p-1}u^{p-1}+O(u^{p})
  \end{displaymath}
  Arguing as above, we find a coordinate $v$ such that $F=v^{-q}$
  around $\infty_2$ and a series
  \begin{displaymath}
    \log{}yv=s_0+s_1v+\dots+s_{q-1}v^{q-1}+O(v^{q})
  \end{displaymath}
  Write $t=(t_1,\dots,t_{p-1})$ and $s=(s_1,\dots,s_{q-1})$. The
  interested reader will check, by series expansion of
  $x=u^{-1}\exp(\sum_{i\geq{}0}t_iu^i)$ and analogously for $y$, the
  following claim:
  \begin{quote}
    \emph{The functions $(\epsilon'=\log\epsilon,t,s,c)$ form a coordinate
    system. The functions $t$, resp. $s$, only depend on $a$, resp.
    $b$.}
  \end{quote}
  We can now show that $\left<,\right>$ has a constant matrix in these
  coordinates. Let us take for example $\cvf{t_i}$. We have
  \begin{equation}
    \label{eq:3}
    \begin{aligned}[b]
      &\frac{1}{x}\cvf[x]{t_i}=u^i,~\frac{1}{y}\cvf[y]{t_i}=-u^i\\
      &\cvf[F]{x}\cvf[x]{t}+\cvf[F]{y}\cvf[y]{t}=\left(x\cvf[F]{x}-y\cvf[F]{y}\right)u^i=Hu^i
    \end{aligned}
  \end{equation}
  As the functions $t$ only depend on $a$, according to (\ref{eq:2})
  we only need to look at the residues at $\infty_2$. Hence
  \begin{equation}
    \label{eq:4}
    \begin{aligned}
      \left<\cvf{t_i},\cvf{t_j}\right>&=\Res_{\infty_2}\frac{(Hu^i)(Hu^j)}{H}\alpha=\Res_{\infty_2}u^{i+j}H\alpha\\
      &=-\Res_{\infty_2}u^{i+j}dF=\Res_{u=0}pu^{i+j-(p+1)}du=p \delta^{i+j}_{p}
    \end{aligned}
  \end{equation}
  A similar calculation, together with the orthogonality relations
  between $a$ and $b$, and hence between $t$ and $s$, deduced from
  (\ref{eq:2}) proves that $\left<,\right>$ is flat.
  
  To finish, we prove the last claim. The Euler vector field
  corresponds to the class of $F$ in $\oo_{\gx,x}$. It is given by
  \begin{equation}
    \label{eq:5}
    E=\left(\frac{1}{p}+\frac{1}{q}\right)\epsilon\cvf{\epsilon}+\sum_{i=1}^{p-1}\frac{p-i}{p}a_{i}\cvf{a_i}+\sum_{i=1}^{q-1}\frac{q-i}{q}b_{i}\cvf{b_i}+c\cvf{c}
  \end{equation}
  Giving weights $1/p+1/q$ to the variable $\epsilon$, $p-i/p$ to
  $a_i$, $q-i/q$ to $b_i$ and $1$ to $c$, we have that a polynomial
  $h(\epsilon,a,b,c)$ is quasi-homogeneous of degree $d$ if and only
  if $\Lie_E(h)=d\cdot{}h$. From (\ref{eq:2}) we see that the entry in
  the position $ij$ of $M_{\infty_2}^{-1}$ (resp.
  $M_{\infty_1}^{-1}$) is (if not constant) quasi-homogeneous of
  degree $(i+j-p)/p$ (resp.  $(i+j-q)/q$). This proves the claim.
\end{proof}
\section{Versal deformations of functions on isolated singularities}
\label{sec:vers-deform-funct}
Given a reduced analytic variety $(X,x)$ and a germ $f\in\mx_{X,x}$,
we will say that $f$ has an isolated singularity if there exists a
representative $f\colon{}U\ra{}S$ onto the complex line $S$ such that
$U\setminus{}\{x\}$ is smooth and $f$ is submersive at any point of
$U\setminus{}\{x\}$. The deformation problem with which we will be
concerned is referred to as \emph{deformations of $X$ over $S$}, that
is, we will consider diagrams
\begin{displaymath}
  \xymatrix{
    (X,x)\ar@{^{(}->}[rr]^{i}\ar[rd]^f\ar[dd]&&(\gx,x)\ar[ld]_F\ar[dd]^{\pi\text{ flat}}\\
    &(S,0)\\
    \{0\}\ar@{^{(}->}[rr]&&(B,0)}
\end{displaymath}
The notions of induced, pull-back or isomorphism of diagrams are
defined in the customary fashion through maps on the base spaces,
keeping the complex line $(S,0)$ fixed. This deformation theory is
sometimes denoted by $\text{Def}(X/S)$ and from a purely algebraic
point of view, it corresponds to the study of the deformations of
$\oo_{X,x}$ as $\oo_{S,0}$ algebra.

As in any deformation theory, we have the powerful theory of the
cotangent cohomology modules at our disposal. Given any holomorphic
map $h\colon{}A\ra{}B$ between analytic spaces, and a $\oo_A$-module
$M$, we will denote by $T^i_{A/B}(M)$ the $i$-th cotangent cohomology
group with coefficients in $M$ (for example \cite{Ill72},
\cite{Ill72bis}. In the absolute case where $B$ reduces to a point, it
is customary to write $T^i_{A}(M)$. If $M$ is just $\oo_A$, then the
notation is further simplified to $T^i_{A/B}$. Another piece of
notation of which we will use is the following. The $0$-th cotangent
cohomology $T^0_{A/B}(M)$ is simply the module of relative vector
fields with coefficients in $M$, that is, $\Theta_{A/B}\otimes_B{}M$.
We will rather use this notation (and those derived, like
$T^0_{B}(\oo_{A})=\Theta(h)$) to be in line with the long established
tradition in Singularity Theory.

Most of the usefulness of the cotangent modules, as for any cohomology
theory, resides in the long exact sequences derived from short exact
sequences of modules. In the particular case of the cotangent
cohomology modules, this is if possible more so. We not only obtain
long exact sequences from short exact sequences of modules, but also
from homomorphism of the base rings (a neat review of the properties
we will use can be found in \cite{BehChr91}). Back to our function
$f\colon{}(X,x)\ra{}(S,0)$, there are two sequences of special
relevance. The first is obtained by considering the problem of
deforming $(X,x)$ alone. If $\pi\colon(\gx,x)\ra{}(B,0)$ is a (flat)
deformation of $(X,x)$, it is the \emph{Zariski-Jacobi long exact
  sequence} associated to the ring homomorphism
$\cc\ra\oo_{X,x}\ra\oo_{\gx,x}$. It begins
\begin{equation}
  \label{eq:6}
  0\ra{}T^0_{\gx/B,x}\ra{}T^0_{\gx,x}\ra{}T^0_{B,0}(\oo_{\gx,x})\ra{}T^1_{\gx/B,x}\ra{}T^1_{\gx,x}\ra{}\dots
\end{equation}
The composite of $\Theta_{B,0}\ra\Theta(\pi)_x$ with the connecting
homomorphism of (\ref{eq:6}) is the \emph{Kodaira-Spencer map} of the
deformation. Its kernel is the submodule of \emph{liftable vector
  fields} and we will denote it by $\gl_{\pi,0}$. In many interesting
cases it coincides with those vector fields tangent to the
discriminant of $\pi$.

If we now consider an extension $F$ of $f$ to the total space
$(\gx,x)$, we can write $\varphi=(\pi,F)$. The second sequence is also
a Zariski-Jacobi sequence, this time corresponding to
$\oo_{B,0}\ra\oo_{S\times{}B,0}\ra\oo_{\gx,x}$:
\begin{equation}
  \label{eq:7}
  0\ra{}T^0_{\gx/S\times{B},x}\ra{}T^0_{\gx/B,x}\ra{}T^0_{S\times{B}/B,0}(\oo_{\gx,x})\ra{}T^1_{\gx/S\times{B},x}\ra{}T^1_{\gx/B,x}\dots
\end{equation}
As before, we will be specially interested in a kernel, this time that
of the map $T^1_{\gx/S\times{B},x}\ra{}T^1_{\gx/B,x}$. We will denote
it by $M_{\varphi,x}$. In fact, this module is readily described in
more familiar terms using the exactness of (\ref{eq:7}). If
$tF\colon\oo_{\gx,x}\ra\Theta(F)_x$ denotes the tangent map of $F$, we
have
\begin{equation}
  \label{eq:8}
  M_{\varphi,x}=\frac{\Theta(F)_x}{tF(\Theta_{\gx/B,x})}
\end{equation}
After all these clarifications, we can state the main lemma of this
section. The proof is so straightforward that it can be safely left to
reader. It neatly separates the problem of finding a versal
deformation of a function on a singular germ into firstly, versally
deforming $(X,x)$ and secondly, versally unfolding $f$.
\begin{lem}
  \label{le:1}
  There is a commutative diagram
  \begin{displaymath}
    \xymatrix{
      0\ar[r]&\gl_{\pi,0}\ar[r]\ar[d]^{-t'F}&\Theta_{B,0}\ar[r]\ar[d]&T^1_{\gx/B,x}\ar@{=}[d]^{\mathbb
        I}\ar[r]&T^1_{\gx,x}\\
      0\ar[r]&M_{\varphi,x}\ar[r]&T^1_{\gx/S\times{B},x}\ar[r]&T^1_{\gx/B,x}\ar[r]&0}
  \end{displaymath}      
  {\noindent}where $t'F$ is defined as follows: for $u\in\gl_{\pi,0}$,
  let $\tilde{u}\in \Theta_{\gx,x}$ be a lift of $u$. Then $t'F(v)$ is
  the class of $tF(\tilde{u})$ in $M_{\varphi,x}$.
\end{lem}
\begin{rem}
  The vertical arrow in the middle is the Kodaira-Spencer map of the
  the map $\varphi$ understood as a deformation of
  $f\colon{}(X,x)\ra(S,0)$.  It is the composite of
  $\Theta_{B,0}\ra\Theta({\varphi})_x$ with the connecting
  homomorphism of Zariski-Jacob sequence derived from
  $\oo_{S,0}\ra\oo_{S\times{}B,0}\ra\oo_{\gx,x}$.
\end{rem}
We deduce the following criterion for versality:
\begin{cor}
  \label{versality-criterion}
  A deformation $\varphi=(F,\pi)$ of $f\colon{}(X,x)\ra(S,0)$ is
  versal if and only if $\pi$ is versal as a deformation of $(X,x)$ and $t'F$ is surjective.
\end{cor}
Versal deformations can be now easily constructed from a versal
deformations $\pi$ of $(X,x)$. We take $f_1,\dots,f_l$ generators of
the vector space $\text{coker}~t'f/\mx_{B,0}(\text{coker}~t'f)$ and
consider the function $F=f+a_1f_1+\dots+a_lf_l$ adding new parameters
$a_1,\dots,a_l$. Requiring that $\pi$ be miniversal and
$f_1,\dots,f_l$ be a basis we will obtain a miniversal deformation.
We will later see examples where this can be explicitly carried out.

\section{Milnor and Tjurina numbers}
\label{sec:miln-tjur-numb}

An important feature of unfoldings of isolated singularities on smooth
spaces is the conservation of the Milnor number. This invariant can be
defined, among other ways, as the length of the Jacobian
$\oo_{\cc^{n+1},0}/(\cvf[f]{x_1},\dots,\cvf[f]{x_{n+1}})$. It is
therefore both the number of non-degenerate critical points of a
generic unfolding and the minimal number of parameters needed to
versally unfold $f$. From this latter point of view, it could also be
called the \emph{Tjurina number} of the deformation problem defined by
right equivalence of functions. 

In our situation, even if the singularity $(X,x)$ is smoothable and we
can speak of non-degenerate critical points of an unfolding, we might
have a different number of those in non-isomorphic Milnor fibres. An
example of this phenomenon is provided by the linear section
$f=x_0+x_1+x_2+x_3+x_4$ on the germ $(X,0)$ of the cone over the
rational normal curve of degree $4$ (\cite{Pin74}). On the other hand,
we do have a well defined Tjurina number as the dimension of the
vector space of first order infinitesimal deformations, namely, the
length $\tau(X/S)$ of $T^1_{X/S,x}$. The next theorem tell us of
conditions in which the Tjurina number indeed coincides with the
number of non-degenerate critical points in every generic deformation.

\begin{prop}
  \label{prop:1}
  Let $\varphi=(F,\pi)\colon(\gx,0)\ra(S\times B,0)$ be a
  $1$-parameter deformation of $f$. Assume that the following
  \emph{extendability condition} is satisfied:
  {\vskip\jot}
  \begin{tabular}[hbpt]{cc}
    \parbox{.05\textwidth}{(*)}&\parbox{.7\textwidth}{\emph{any vector field tangent to the fibres of $f$ can be extended to a vector field tangent to the fibres of $\varphi$,}}
  \end{tabular}
  {\vskip\jot}
  \noindent{}then both $T^1_{\gx/S\times{}B,x}$ and $M_{\varphi,x}$ are free
  $\oo_{B,0}$-modules. Moreover, if $T^2_{X,x}=0$ and the generic
  fibre of $\pi$ is smooth, their ranks coincide.
\end{prop}
\begin{proof}   
  Let $y$ be a parameter in $(B,0)$. The exact sequence
  \[0\ra{}\oo_{\gx,x}\xra{\cdot~y}\oo_{\gx,x}\ra\oo_{X,x}\ra{}0\]
  {\noindent}induces a long exact sequence:
  \begin{align}
    \label{eq:9}
    \begin{split}
      0&\ra{}\Theta_{\gx/S\times{B},x}\xra{\cdot{y}}\Theta_{\gx/S\times{B},x}\lra{}\Theta_{X/S,x}\\
      &\ra{}T^1_{\gx/S\times{}B,x}\xra{\cdot{}y}T^1_{\gx/S\times{}B,x}\ra{}T^1_{X/S,x}\ra\dots
      \end{split}
  \end{align}  
  The condition (*) implies that map
  $\Theta_{\gx/S\times{}B,x}\ra{\Theta_{X/S,x}}$ is surjective and
  hence $T^1_{\gx/S\times{B},x}\xra{\cdot{y}}T^1_{\gx/S\times{B},x}$
  injective. Therefore $T^1_{\gx/S\times{}B,x}$ and $M_{\varphi,x}$
  are flat over $\cc\{y\}$ and hence free.
  
  For the second statement, we first show that the condition
  $T^2_{X,x}=0$ also implies $T^2_{X/S,x}=0$. Associated to
  $\cc\rightarrow\oo_{S,0}\ra\oo_{X,x}$ we have a long exact sequence:
  \[\ldots\ra{}T^i_{X/S,x}\ra{T^i_{X,x}}\ra{T^i_{S}(\oo_{X,x})}\ra{T^{i+1}_{X/S,x}}\ra\ldots\]
  As $(S,0)$ is smooth, $T^i_{S}(\oo_{X,x})=0$ for $i\geq{1}$, so that
  $T^i_{X/S,x}=T^i_{X,x}$ for $i\geq{2}$. Finally if the generic fibre
  of $g$ is a smooth, then $T^2_{\gx/S\times{B},x}$ is annihilated by
  a power of the maximal ideal $\mx_{B,0}$, and hence it is Artinian.
  The exact sequence (\ref{eq:9}) then contains the short exact
  sequence:
  \begin{displaymath}
    0\ra{}T^1_{\gx/S\times{B},x}\xra{\cdot{y}}T^1_{\gx/S\times{B},x}\ra{}T^1_{X/S,x}\ra{0}
  \end{displaymath}
  It follows that $\rank{T^1_{\gx/S\times{}B,x}}=\dim_\cc~T^1_{X/S,x}$.
  To see that this is also the rank of $M_{\varphi,x}$ we write one
  more exact sequence:
  \begin{displaymath}
    0\ra{}M_{\varphi,x}\ra{T^1_{\gx/S\times{}B,x}}\ra{T^1_{\gx/B,x}}\ra{0}
  \end{displaymath}
  {\noindent}and notice that $T^1_{\gx/B,x}$ is supported at $0$.
\end{proof}
From now on, we restrict ourselves to situations in which all the
conditions of the above theorem are satisfied, namely, functions on
\emph{smoothable and unobstructed singularities} for which the
condition (*) holds \emph{for any $1$-parameter deformation}. We will
now show that this family of functions includes some interesting
examples. First, note that it follows from the above proposition that
not only $\tau(X/S)$ coincides with the number of Morse critical
points in the generic deformation, but also that for the miniversal
deformation of $f$, the map
\begin{equation}
  \label{eq:10}
  t'F\colon{}\gl_{\pi,0}\lra{}M_{\varphi,x}
\end{equation}
\noindent{}extends to an isomorphism of \emph{free sheaves}. In particular,
the sheaf of liftable vector fields is necessarily free.
We will now have a close look at two situations for which we can prove
the extendability condition: the case described in the previous remark
and that of isolated complete intersection singularities. But let us
first introduce a piece of notation. The module $M_{\varphi,x}$ is not
independent of the given deformation $\varphi$. Even its length is not
a well-defined invariant of the function $f$. To avoid such a
dependence, we consider the miniversal deformation of $(X,x)$ alone,
say $\pi\colon(\gx,0)\ra{}(B,0)$ and take any extension $F$ to the
total space. We define
\begin{displaymath}
  M_f=\frac{M_{\varphi,x}}{\mx_{B,0}M_{\varphi,x}}
\end{displaymath}
Note that this module is well-defined as any two extensions of $f$
differ by an element of the maximal ideal $\mx_{B,0}$.

The reason to introduce this module is that if the conditions of
proposition (\ref{prop:1}) are fulfilled, its length will be equal to
$\tau(X/S)$.

\section{Functions on space curves and complete intersections}
\label{sec:funct-space-curv}

As remarked before, the case of functions on smoothable and
unobstructed curves falls trivially into our area of interests. The
equality between the Milnor number and Tjurina number for these is the
main result of \cite{MonStr01}.

The authors of \cite{MonStr01} define the Milnor number of a function
on a space curve in terms of the \emph{dualising module}
$\omega_{X,x}$.  Using the class map (\cite{AngLej89}), or
equivalently, Rosenlicht's description of $\omega_{X,x}$ as
certain meromorphic forms (\cite{BucGre80}), they define the Milnor
number of $f$ as
\begin{displaymath}
  \mu_f=\dim_\cc\frac{\omega_{X,x}}{\oo_{X,x}df}
\end{displaymath}
They also show the following interesting formula: the class map
$cl\colon\Omega_{X,x}\ra\omega_{X,x}$ can be dualised to obtain a
submodule $\omega_{X,x}^*$ of $\Theta_{X,x}$. Then $\mu_f$ is also the
length of $\Theta(f)_x/tf(\omega^*_{X,x})$.
\begin{prop}
  For a function $f$ on a space curve,
  \begin{displaymath}
    M_{f}=\frac{\Theta(f)_x}{tf(\omega^*_{X,x})}
  \end{displaymath}
\end{prop}
\begin{proof}
  A space curve is a Cohen-Macaulay variety of codimension $2$ and as
  such is defined by the maximal minors $\Delta_i$ of a
  $(m\times{}m-1)$-matrix $M$ with coefficients in $\oo_{\cc^3,x}$.
  There deformations are well understood (\cite{Sch77}), they are also
  defined by the maximal minors $\tilde{\Delta}_i$ of a perturbation
  $\tilde{M}$ of $M$.
  
  An identical calculation to that of \cite{MonStr01}, but using the
  relative class map for the miniversal family instead of that of
  $(X,x)$, shows that its dual in $\Theta_{\gx/B,x}$ is generated by
  the vector fields
  \begin{displaymath}
    \begin{vmatrix}
      \displaystyle\cvf{x_1}&\displaystyle\cvf{x_{2}}&\displaystyle\cvf{x_{3}}\\
      \vspace{-3\jot}\\
      \displaystyle\cvf[\tilde{\Delta}_i]{x_1}&\displaystyle\cvf[\tilde{\Delta}_i]{x_{2}}&\displaystyle\cvf[\tilde{\Delta}_i]{x_{3}}\\
      \vspace{-3\jot}\\
      \displaystyle\cvf[\tilde{\Delta}_j]{x_1}&\displaystyle\cvf[\tilde{\Delta}_j]{x_{2}}&\displaystyle\cvf[\tilde{\Delta}_j]{x_{3}}\\
    \end{vmatrix}
    ~,~1\leq{}i<j\leq{}l
  \end{displaymath}
  On the other hand, the relative class map
  $cl_{\gx/B,x}\colon\Omega_{\gx/B,x}\ra\omega_{\gx/B,x}$ (or rather,
  a representative of it) is an isomorphism whenever the fibre is
  smooth.  As the generic fibre is indeed smooth, the set where it
  fails to be bijective is of codimension at least $2$. Hence its dual
  is an isomorphism everywhere and if $F$ is any extension of $f$ to
  $(\gx,x)$ we have
  \begin{displaymath}
    M_f=\Theta(F)_x/tF(\Theta_{\gx/B,x})+\mx_{B,0}\Theta(F)_x=\frac{\Theta(f)_x}{tf(\omega^*_{X,x})}
  \end{displaymath}
\end{proof}
\begin{exmp}
  We can use the above calculation to compute versal deformations of
  functions on space curves. For example, the union of the three
  coordinate axis in $(\cc^3,0)$ is defined by the $2\times{}2$-minors
  of
  $M=\left(\begin{smallmatrix}x&y&0\\0&y&z\end{smallmatrix}\right)$.
  The miniversal deformation of a function $f=x^p+y^q+z^r$ is
  therefore obtained by considering the miniversal deformation of the
  curve together with the unfolding
  $F=x^p+\sum_{i=1}^{p-1}a_{i}x^{p-i}+y^p+\sum_{i=1}^{q-1}b_{i}y^{q-i}+z^r+\sum_{i=1}^{r-1}c_{i}z^{r-i}+d$.
  In \cite{Gor00}, where simple functions on curves are classified, this
  singularity is referred to as $C_{p,q,r}$.
\end{exmp}
We now go on to study the case of functions on complete intersections.
Let $f\colon(X,x)\ra(S,0)$ be a germ with an isolated singularity on a
$n$-dimensional complete intersection.  Let $g_1,\dots,g_k$ be elements
defining the ideal of $(X,x)$ in $(\cc^{n+k},x)$.

A submodule of $\Theta_{X,x}$ whose members are clearly tangent to all
the fibres of $f$ is generated by the maximal minors of the matrix
\begin{equation}
  \label{eq:11}
  \begin{pmatrix}
      \displaystyle\cvf{x_1}&\dots&\displaystyle\cvf{x_{n+k}}\\
      \vspace{-3\jot}\\
      \displaystyle\cvf[f]{x_1}&\dots&\displaystyle\cvf[f]{x_{N}}\\
      \vspace{-3\jot}\\
      \displaystyle\cvf[g_1]{x_1}&\dots&\displaystyle\cvf[g_1]{x_{n+k}}\\
      \hdotsfor{3}\\
      \vspace{-3\jot}\\
      \displaystyle\cvf[g_k]{x_1}&\dots&\displaystyle\cvf[g_k]{x_{n+k}}
    \end{pmatrix}
\end{equation}
\begin{lem}\label{le:3}
  The vector fields in (\ref{eq:11}) generate $\Theta_{X/S}$.
\end{lem}
\begin{proof}
  Let $\varphi=(f,g_1,\dots,g_k)$. The module $\Theta_{X/S,x}$ is the kernel
  of
  \begin{equation}
    \label{eq:12}
    \Theta_{\cc^{n+k},x}\otimes\oo_{X,x}\xra{t\varphi\otimes{}\mathbb{I}}\Theta_{\cc^{1+k},0}\otimes\oo_{X,x}
  \end{equation}
  As $f$ has an isolated singularity and $(X,x)$ is a Cohen-Macaulay,
  the depth of the ideal in $\oo_{X,x}$ generated by the maximal
  minors of (\ref{eq:12}) is $n$, i.e, the greatest possible.  It follows
  that the Eagon-Northcott complex is exact (\cite{EagNor62}) and the
  kernel is generated by the above vector fields.
\end{proof}
\begin{cor}
  For germ $f\colon(X,x)\ra(S,0)$ with an isolated singularity on a
  complete intersection, $\tau(X/S)$ coincides with the number of
  non-degenerate critical points of a generic deformation. If $n$
  denotes the dimension of $(X,x)$, then
  \begin{equation}
    \label{eq:13}
    M_f\simeq\frac{\omega_{X,x}}{df\wedge\Omega_{X,x}^{n-1}}
  \end{equation}
\end{cor}
\begin{proof}
  In view of lemma (\ref{le:3}), it is clear that the extendability
  condition holds. But lemma (\ref{le:3}) is also telling us which
  vector fields are tangent to all the fibers of a deformation of a
  complete intersection. Simply take $(X,x)$ to be the ambient space
  $(\cc^{n+k},x)$ and change $(S,0)$ by $(\cc^{k},0)$. We see that
  $\Theta_{\gx/B,x}$ is generated by the maximal minors of (\ref{eq:11})
  with the row involving $f$ deleted. The equality (\ref{eq:13}) is now
  evident by differentating $f$ with respect to this set of generators
  of $\Theta_{\gx/B,x}$.
\end{proof}
\begin{rem}
  The equality between Tjurina and Milnor numbers for functions on
  complete intersections is proven with unrelated methods in
  \cite{Gor95}.
\end{rem}
\begin{rem}
  Using (\ref{eq:13}) and a well-known result \cite{LeDung78}, we can
  interpret the rank of $M_f$, and hence $\tau(X/S)$, as the rank of
  certain vanishing homology, namely $H_{n}(X_b,Y_s)$ for Milnor
  fibers of $(X,x)$ and $f$.
\end{rem}

\section{Multiplication on the sheaf of liftable vector fields}
\label{sec:mult-sheaf-lift}

Whenever the map $t'F$ of (\ref{eq:10}) extends to an isomorphism of
sheaves, we can use it to define a multiplication on $\gl_{\pi}$ by
pulling-back the algebra structure on $M_{\varphi}$. If $(X,x)$ is
also smoothable then this defines a multiplication on the tangent
bundle of the complement $B\setminus{}\Delta$ of the discriminant of
the fibration. We begin recalling the definition of $F$-manifold from
\cite{Her02}~Ch.~1.
\begin{defn}
  A complex manifold with an associative and commutative multiplication
$\star$ on the tangent bundle is called an $F$-manifold if: 
\begin{enumerate}
  \item \emph{(unity)} there exists a \emph{global} vector field $e$
    such that $e\star{}u=u$ for any $u\in\Theta_M$ and,
    \item\emph{(integrability)}
      $\Lie_{u\star{}v}(\star)=u\star{}\Lie_v(\star)+\Lie_u(\star)\star{}v$ for any $u,v\in\Theta_M$.
\end{enumerate}
An \emph{Euler vector field} $E$ (of weight $1$) for $M$ is defined by
the condition
\begin{displaymath}
  \Lie_E(\star)=\star
\end{displaymath}
\end{defn}
The main consequence of this definition is the integrability of
multiplicative subbundles of $TM$, namely, if in a neigbourhood $U$ of
a point $p\in{}M$ we can decompose $TU$ as a sum of unitary
subalgebras $A\oplus{}B$ such that $A\star{}B$, then $A$ and $B$ are
integrable.

By choosing good representatives in the sense of \cite{Loo84} for all the
germs involved, we have:
\begin{prop}
  \label{prop:2}
  The map $t'F$ endows the sheaf of liftable vector fields
  $\gl_{\pi}$ with the structure of commutative and associative
  $\oo_B$-algebra $\star$ such that, for any $u,v\in\gl_{\pi}$:
  \begin{equation}
    \label{eq:14}
    \Lie_{u\star{}v}(\star)=\Lie_u(\star)\star{}v+u\star\Lie_v(\star)
  \end{equation}
  The class of $F$ in $\pi_*M_{\varphi}$ corresponds to an Euler vector field of weight $1$.
\end{prop}
\begin{proof}
  It is enough to show (\ref{eq:14}) off $\Delta$. Let
  $\mu=\rank{}\pi_*M_\varphi$.  For a generic point
  $b\in{}B\setminus{}\Delta$, the function $F$ has $\mu$ quadratic
  singularities on the smooth fiber $\pi^{-1}(b)$. Hence
  $\pi_*M_{\varphi}$ decomposes into $\mu$ 1-dimensional unitary
  subalgebras. In a neigbourhood $U\subset{}S\setminus{}\Delta$ of
  such a point, the integrability condition is equivalent to the image
  $L$ of the map
  \begin{equation}
    \label{eq:15}
    \text{supp}~M_\varphi\ni{}x\mapsto{}d_xF\in{}T_{\pi(x)}^*B
  \end{equation}
  \noindent{}being a Lagrangian subvariety of $T^*B$ (see
  \cite{Her02},~Th.~3.2). If $\alpha$ denotes the canonical $1$-form on
  $T^*B$ and $p\colon{}T^*B\rightarrow{}B$ the projection, it is easy
  to check that the diagram
  \begin{displaymath}
    \xymatrix@-1pc{\text{supp}~M_\varphi\ar[rr]&&p_*\mathcal{O}_L\\
      &\Theta_B\ar[lu]\ar[ru]}
  \end{displaymath}
  \noindent{}is commutative. The homomorphism on the right hand side
  is given by evaluation, so that it can also be expressed as
  $\alpha(\tilde{u})$ where $\tilde{u}$ is a lift of $u\in\Theta_B$ to
  $\Theta_{T^*B}$. Hence $\alpha_L$ is the relative differential of
  $F$ when thought of as a map on $L$ via the identification
  (\ref{eq:15}). It follows that $\alpha_L$ is the
  exact and hence closed, so that $L$  is Lagrangian.\\
  The statement about the Euler vector is an easy calculation that we
  leave to the reader (see \cite{Her02}, Th. 3.3.).
\end{proof}
The above proposition establishes the structure of $F$-manifold at
least off $\Delta$. In the case where $\gl_\pi$ coincides with the
sheaf of tangent vector fields to $\Delta$, denoted by
$\Theta(-\log\Delta)$, we can in fact define the $F$-manifold structure
on each of the stratum of the logarithmic stratification induced by
$\Theta(-\log\Delta)$ (\cite{KSai80}). First we need a lemma:
\begin{lem}
  \label{le:4}
  For any ideal sheaf $I\subset{}\oo_{B}$, the kernel of the map
  \begin{displaymath}
    \gl_{\pi}/I\gl_{\pi}\lra{}\Theta_{B}/I\Theta_{B}
  \end{displaymath}
  \noindent{}is identified by $t'F$ with an ideal of
  $\pi_*M_{\varphi}/I\pi_*M_{\varphi}$.
\end{lem}
\begin{proof}
  Choosing a Stein representative of $\pi$, we can interpret the
  diagram in lemma (\ref{le:1}) as a morphism of free resolutions of
  $\pi_*T^1_{\gx/B}$ where the row at the bottom is actually a complex of
  $g_*\oo_{\gx,x}$-modules. Hence
  \begin{displaymath}
    \text{ker}~(\gl_{\pi}/I\gl_{\pi}\lra{}\Theta_B/I\Theta_B)=\text{Tor}_1^{\oo_{B}}(\pi_*T^1_{\gx/B},\oo_{B}/I)
  \end{displaymath}
  \noindent{}is mapped by $t'F$ to the kernel of
  \begin{displaymath}
    \pi_*M_{\varphi}/I\pi_*M_{\varphi}\ra{}\pi_*T^1_{\gx/S\times{}B}/I\pi_*T^1_{\gx/S\times{}B}
  \end{displaymath}
  \noindent{}and therefore it is an ideal.
\end{proof}
\begin{thm}
  \label{thm:1}
  If $\gl_{\pi}=\Theta(-\log\Delta)$, then each stratum of the
  logarithmic stratification has the structure of $F$-manifold with an
  Euler vector field of weight $1$.
\end{thm}
\begin{proof}
  Let $b\in{}B$ and let $S_b$ be the stratum in which $b$ lies. Let
  $V$ be a open neighbourhood of $b$ in which $S_b\cap{}V$ is an
  analytic subset of $V$ defined by the ideal $I_{S_b}$. The sheaf
  $\Theta_{S_b\cap{}V}$ can be identified with
  \begin{displaymath}
    \frac{\text{im}~(\gl_{\pi}|_{V}\ra\Theta_B|_{V})}{I_{S_b}\text{im}~(\gl_{\pi}|_{V}\ra\Theta_B|_{V})}
  \end{displaymath}
  Let $\mathcal{K}$ denote the sheaf
  $\text{Tor}_1^{\oo_{B}}(\pi_*T^1_{\gx/B},\oo_{B}/I_{S_b})$. The map
  $t'F$ descends to the above quotient and it yields an isomorphism of
  $\oo_{S_b\cap{}V}$-modules
  \begin{displaymath}
    \Theta_{S_b\cap{}V}\xra{~\simeq~}\frac{\pi_*M_{\varphi}|_{V}}{I_{S_b}\pi_*M_{\varphi}|_{V}+t'F(\mathcal{K}|_{V})}
  \end{displaymath}
  According to the previous lemma, the right hand side is a
  $\oo_B$-algebra. The above isomorphism defines the multiplication on
  the tangent bundle of the stratum $S_b$. From propostion
  (\ref{prop:2}) it follows that it is an $F$-manifold with
  Euler vector field of weight $1$ given by the class of $F$ in the
  corresponding algebra.
\end{proof}
\begin{rem}
  If the stratum $S_b$ is a massive $F$-manifold, that is, there exist
  coordinates $u_1,\dots,u_l$ such that
  \begin{displaymath}
    \cvf{u_i}\star\cvf{u_j}=\delta_{ij}\cvf{u_i}~\text{for all $i,j$}
  \end{displaymath}
  \noindent{} then the critical values of $F$ are generically local 
  coordinates on $S_b$. In particular, this always holds on the
  stratum $B-\Delta_g$. In the case of space curves, that the critical
  values of $F$ off the bifurcation diagram are local coordinates is
  shown for simple functions in \cite{Gor00bis} who also conjectured
  the analogous result for non-simple functions.
\end{rem}


\bibliographystyle{amsplain} \bibliography{Bibliography}

\def\cprime{$'$}
\providecommand{\bysame}{\leavevmode\hbox to3em{\hrulefill}\thinspace}
\providecommand{\MR}{\relax\ifhmode\unskip\space\fi MR }
\providecommand{\MRhref}[2]{%
  \href{http://www.ams.org/mathscinet-getitem?mr=#1}{#2}
}
\providecommand{\href}[2]{#2}
\begin{thebibliography}{10}

\bibitem{AngLej89}
B.~Ang{\'e}niol and M.~Lejeune-Jalabert, \emph{Calcul diff\'erentiel et classes
  caract\'eristiques en g\'eom\'etrie alg\'ebrique}, Travaux en Cours [Works in
  Progress], vol.~38, Hermann, Paris, 1989, With an English summary.
  \MR{90h:14004}

\bibitem{Bar00}
Serguei Barannikov, \emph{Semi-infinite {H}odge structures and mirror symmetry
  for projective spaces}, AG/0108148.

\bibitem{BehChr91}
Kurt Behnke and Jan~Arthur Christophersen, \emph{Hypersurface sections and
  obstructions (rational surface singularities)}, Compositio Math. \textbf{77}
  (1991), no.~3, 233--268. \MR{MR1092769 (92c:14002)}

\bibitem{BucGre80}
Ragnar-Olaf Buchweitz and Gert-Martin Greuel, \emph{The {M}ilnor number and
  deformations of complex curve singularities}, Invent. Math. \textbf{58}
  (1980), no.~3, 241--281. \MR{81j:14007}

\bibitem{Dub96}
Boris Dubrovin, \emph{Geometry of {$2$}{D} topological field theories},
  Integrable systems and quantum groups (Montecatini Terme, 1993), Lecture
  Notes in Math., vol. 1620, Springer, Berlin, 1996, pp.~120--348.
  \MR{97d:58038}

\bibitem{EagNor62}
J.~A. Eagon and D.~G. Northcott, \emph{Ideals defined by matrices and a certain
  complex associated with them.}, Proc. Roy. Soc. Ser. A \textbf{269} (1962),
  188--204. \MR{26 \#161}

\bibitem{Gor95}
V.~V. Goryunov, \emph{Singularities of projections}, Singularity theory
  (Trieste, 1991), World Sci. Publishing, River Edge, NJ, 1995, pp.~229--247.
  \MR{96m:58018}

\bibitem{Gor00bis}
\bysame, \emph{Functions on space curves}, J. London Math. Soc. (2) \textbf{61}
  (2000), no.~3, 807--822. \MR{2002f:58068}

\bibitem{Gor00}
\bysame, \emph{Simple functions on space curves}, Funktsional. Anal. i
  Prilozhen. \textbf{34} (2000), no.~2, 63--67. \MR{2001g:58065}

\bibitem{HerMan99}
C.~Hertling and Yu. Manin, \emph{Weak {F}robenius manifolds}, Internat. Math.
  Res. Notices (1999), no.~6, 277--286. \MR{2000j:53117}

\bibitem{Her02}
Claus Hertling, \emph{Frobenius manifolds and moduli spaces for singularities},
  Cambridge Tracts in Mathematics, vol. 151, Cambridge University Press,
  Cambridge, 2002. \MR{1 924 259}

\bibitem{Ill72bis}
Luc Illusie, \emph{Complexe cotangent et d\'eformations. {I}}, Springer-Verlag,
  Berlin, 1971, Lecture Notes in Mathematics, Vol. 239. \MR{58 \#10886a}

\bibitem{Ill72}
\bysame, \emph{Complexe cotangent et d\'eformations. {II}}, Springer-Verlag,
  Berlin, 1972, Lecture Notes in Mathematics, Vol. 283. \MR{58 \#10886b}

\bibitem{Loo84}
E.~J.~N. Looijenga, \emph{Isolated singular points on complete intersections},
  London Mathematical Society Lecture Note Series, vol.~77, Cambridge
  University Press, Cambridge, 1984. \MR{86a:32021}

\bibitem{MonStr01}
David Mond and Duco van Straten, \emph{Milnor number equals {T}jurina number
  for functions on space curves}, J. London Math. Soc. (2) \textbf{63} (2001),
  no.~1, 177--187. \MR{2002e:32040}

\bibitem{Pin74}
Henry~C. Pinkham, \emph{Deformations of algebraic varieties with {$G\sb{m}$}
  action}, Soci\'et\'e Math\'ematique de France, Paris, 1974, Ast\'erisque, No.
  20. \MR{51 \#12847}

\bibitem{Ros98}
Markus Rosellen, \emph{{H}urwitz spaces and {F}robenius manifolds}, preprint
  Max-Planck-Institut für Mathematik, 97--99.

\bibitem{KSai80}
Kyoji Saito, \emph{Theory of logarithmic differential forms and logarithmic
  vector fields}, J. Fac. Sci. Univ. Tokyo Sect. IA Math. \textbf{27} (1980),
  no.~2, 265--291. \MR{83h:32023}

\bibitem{Sch77}
Mary Schaps, \emph{Deformations of {C}ohen-{M}acaulay schemes of codimension
  {$2$} and non-singular deformations of space curves}, Amer. J. Math.
  \textbf{99} (1977), no.~4, 669--685. \MR{58 \#10918}

\bibitem{LeDung78}
L{\^e}~D ung Tr{\'a}ng, \emph{The geometry of the monodromy theorem}, C. P.
  Ramanujam---a tribute, Tata Inst. Fund. Res. Studies in Math., vol.~8,
  Springer, Berlin, 1978, pp.~157--173. \MR{80i:32031}

\end{thebibliography}

\end{document}